\documentclass{amsart}
\usepackage{amsmath}
\usepackage{amsthm}
\usepackage{amssymb}
\usepackage{graphicx}
\usepackage{todonotes}
\usepackage{enumitem}
\usepackage{soul}

\newtheorem{assumption}{Assumption}
\newtheorem{theorem}{Theorem}[section]

\newtheorem{lemma}[theorem]{Lemma}

\newtheorem{definition}{Definition}
\newtheorem{remark}{Remark}

\newenvironment{vf}{\left\{\begin{array}{rcl}}{\end{array}\right.}
\DeclareMathOperator{\divergenceOperator}{div}

\newcommand{\includegraph}[2][]{\ifnum\pdfoutput=0\includegraphics[#1]{#2.eps}\else\includegraphics[#1]{#2.pdf}\fi}

\newcommand{\R}{\mathbb{R}}
\newcommand{\T}{\mathbb{T}}
\renewcommand{\S}{\mathbb{S}}

\newcommand{\rev}[1]{{\textcolor{black}{#1}}}

\title{Slow-fast torus knots}
\author{Renato Huzak}
\address[R. Huzak]{Hasselt University, Campus Diepenbeek, Agoralaan Gebouw D, 3590 Diepenbeek, Belgium}
\email{renato.huzak@uhasselt.be}
\author{Hildeberto Jardón-Kojakhmetov}
\address[H. Jardón-Kojakhmetov]{University of Groningen,
Faculty of Science and Engineering,
Dynamical Systems, Geometry \& Mathematical Physics — Bernoulli Institute,
Nijenborgh 9,
9747 AG, Groningen,
The Netherlands}
\email{h.jardon.kojakhmetov@rug.nl}
\date{}

\begin{document}



\maketitle
\begin{abstract}
The goal of this paper is to study global dynamics of $C^\infty$-smooth slow-fast systems on the $2$-torus of class $C^\infty$ using geometric singular perturbation theory and the notion of slow divergence integral. Given any $m\in\mathbb{N}$ and two relatively prime integers $k$ and $l$, we show that there exists a slow-fast system $Y_{\epsilon}$ on the $2$-torus that has a $2m$-link of type $(k,l)$, i.e. a (disjoint finite) union of $2m$ slow-fast limit cycles each of $(k,l)$-torus knot type, for all small $\epsilon>0$. The $(k,l)$-torus knot turns around the $2$-torus $k$ times meridionally and $l$ times longitudinally. There are exactly $m$ repelling 
limit cycles and $m$ attracting 
limit cycles. \rev{Our analysis: a) proves the case of normally hyperbolic singular knots, and b) provides sufficient evidence to conjecture a similar result in some cases where the singular knots have regular nilpotent contact with the fast foliation.}

\end{abstract}




\section{Introduction}\label{Section-Introduction}
Singularly perturbed systems on the $2$-torus $\mathbb{T}^2$ have been studied in \cite{GuIl} and \cite{Schurov4,Schurov2,Schurov1,Schurov3}. In \cite{GuIl} the authors constructed a slow-fast system on $\mathbb{T}^2$ (depending only on a singular parameter $\epsilon$) with the following property: there is a sequence of $\epsilon$-intervals accumulating at $0$ such that the system has exactly $2$ limit cycles (one is a stable canard and the other one is an unstable canard) for each $\epsilon>0$ from these intervals. \rev{In \cite{GuIl}, a} limit cycle is called a canard
limit cycle if it contains a part passing near repelling portions of the critical \rev{manifold}. 
\rev{A more classical definition of canard is as follows \cite{de2021canard}: a canard limit cycle is a closed orbit that passes near \emph{both} attracting and repelling portions of the critical manifold. This is the definition we adopt through this article. So, in such context, the repelling limit cycles obtained below, are not canards. We further notice that our framework is essentially different from the torus canards addressed in e.g. \cite{vo2017generic} and related works.}

In \cite{Schurov4,Schurov2,Schurov1,Schurov3}, I. V. Schurov generalized the results of \cite{GuIl} (the main focus has been directed towards the existence of (attracting) canard cycles, in the sense of \cite{GuIl}, on $\mathbb{T}^2$ under more general conditions). The main tool is the Poincar\'{e} map from $\mathbb{S}^1$ to itself. If the rotation number of the Poincar\'{e} map is an integer and the critical manifold is connected, then the number of canard limit cycles is bounded by the number of fold points of the critical manifold (see \cite{Schurov2}).
\smallskip

The papers mentioned above mainly deal with ``unknotted" limit cycles on the $2$-torus (the limit cycles make one pass along the slow direction, i.e. the case of integer rotation number) and the connected critical \rev{manifold} with jump contact points does not turn around the $2$-torus along \rev{both} the slow direction and the fast direction. An exception is \cite{Schurov4} where a result was proved for non-integer rotation number (more precisely, canard limit cycles -in the sense of \cite{GuIl}- make two passes along the slow direction). \textit{The main purpose of our paper is to show the existence of slow-fast systems, with one \rev{small} parameter $\epsilon$, with an arbitrary finite number of repelling \rev{limit} cycles on $\mathbb{T}^2$ that make $k$ passes along the slow direction and $l$ passes along the fast direction ($k$ and $l$ are relatively prime)}. Such limit cycles are called $(k,l)$-torus knots and occur for each small $\epsilon>0$ (for more details about the torus knots see the rest of this section and Appendix \ref{sec:app1}). They are generated by an unconnected normally hyperbolic critical curve (each component is a $(k,l)$-torus knot). In order to prove this and to find the global dynamics on $\mathbb{T}^2$, we will use a fixed point theorem for segments, the notion of slow divergence integral defined along the critical \rev{manifold} and a generalization of \rev{the} Poincar\'{e}-Bendixson theorem to $\mathbb{T}^2$ due to Schwartz (see \cite{Schw}). In our slow-fast setting, it is more convenient to not use the Poincar\'{e} map (see Remark \ref{remark-Poincare}). 

\smallskip

\begin{figure}[htb]
    \begin{center}
        \ifnum\pdfoutput=0\includegraphics[bb=0 0 1302 839,width=11cm,height=5cm]{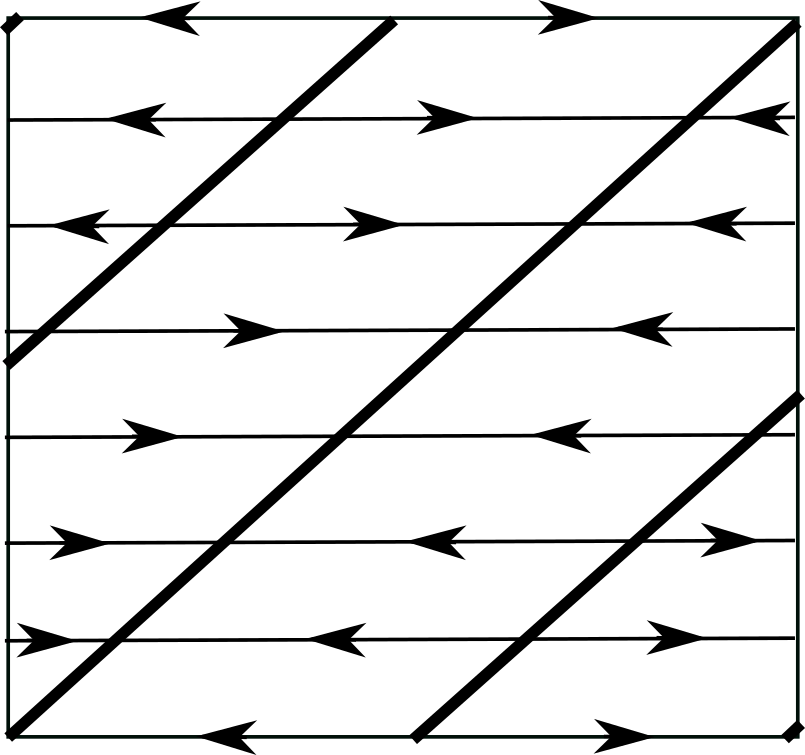}\else\includegraphics[width=4.5cm,height=4.5cm]{torusexample.png}\fi
 {\footnotesize
            \put(-134,-2){$0$} \put(-53,27){$C_+$}  \put(-87,57){$C_-$}  \put(-110,96){$C_+$} \put(-1,-2){$2\pi$} \put(-134,127){$2\pi$}}
 \end{center}
    \caption{Dynamics of $X_0$.}
    \label{figure:torusexample}
\end{figure}

To start fixing ideas, let us consider a slow-fast system
\begin{equation}
\label{equ100}
X_\epsilon: \begin{vf}
        \dot{x} &=&\sin{(y-x)}   \\
        \dot{y} &=&\epsilon
    \end{vf}
\end{equation}
where $(x,y)\in \mathbb{T}^2=\mathbb{R}^2\setminus (2\pi\mathbb{Z}^2)$,  $\epsilon\ge 0$ is a (small) singular parameter\rev{, and the over-dot denotes differentiation with respect to the fast time scale $t$}. The vector field $X_{\epsilon}$ is $2\pi$-periodic in both variables and we restrict our attention to the dynamics of $X_\epsilon$, with $\epsilon\ge 0$, on the two-dimensional torus $\mathbb{T}^2$ (we keep $x,y$ in $[0,2\pi[$ and glue together the opposite segments $x=0$ and $x=2\pi$, and $y=0$ and $y=2\pi$). In the limit $\epsilon=0$, system \eqref{equ100} has horizontal fast orbits and two disjoint simple closed curves of singularities given by $C_-=\{y=x\}$ and $C_+=\{y=x+\pi\}$. These two closed curves pass through $\mathbb{T}^2$ horizontally and vertically only once (see Figure \ref{figure:torusexample}). All the singularities on $C_-$ (resp. $C_+$) are normally attracting (resp. normally repelling). When $\epsilon>0$, these singularities disappear and the dynamics near $C_\pm$ are given by the reduced \rev{slow} flow \rev{$y'=1$, where now the prime denotes differentiation with respect to the slow-time $s=\epsilon t$.}

We are interested in the dynamics of the regular system $X_{\epsilon}$ on $\mathbb{T}^2$ for each small and positive parameter $\epsilon$. Any orbit $\mathcal{O}_{\epsilon}$ of $X_{\epsilon}$ with initial point located away from $C_\pm$ is attracted to $C_-$ and stays close to $C_-$ forever. Such \rev{an} orbit $\mathcal{O}_\epsilon$ cannot therefore be closed. It is clear now that closed orbits of $X_\epsilon$ may appear only in a tubular neighborhood of $C_-$ or $C_+$. Notice that $C_-$ and $C_+$ are limit periodic sets at \rev{the} level $\epsilon=0$ without fast segments ($C_-$ consists of one attracting slow part while $C_+$ has one repelling slow part). 
\smallskip

We show that $C_-$ (resp. $C_+$) generates precisely one limit cycle which is hyperbolic and attracting (resp. repelling), for each small and positive $\epsilon$. This result will be true in an $\epsilon$-uniform neighborhood of $C_-$ (resp. $C_+$). It means that the neighborhood does not shrink to $C_\pm$ as $\epsilon\to 0$. The $\omega$-limit set (resp. the $\alpha$-limit set) of all other orbits (different from the two limit cycles) is the attracting (resp. repelling) limit cycle Hausdorff close to $C_-$ (resp. $C_+$). 
\smallskip

\begin{figure}[htb]
    \begin{center}
        \ifnum\pdfoutput=0\includegraphics[bb=0 0 1302 839,width=11cm,height=5cm]{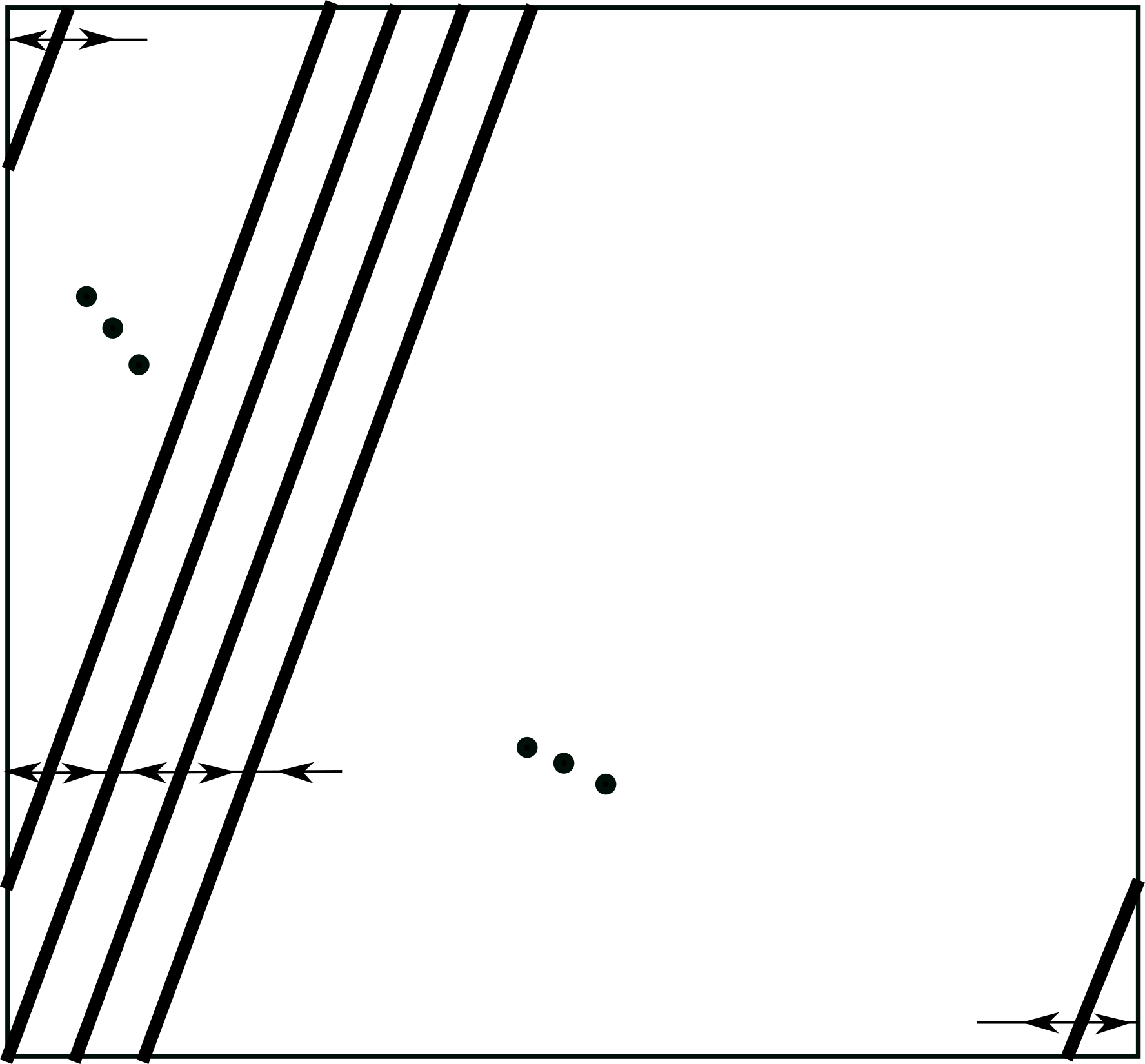}\else\includegraphics[width=5.5cm,height=5.5cm]{torus-general.png}\fi
 {\footnotesize
            \put(-162,-7){$0$}\put(-169,18){$\frac{\pi}{ml}$} \put(-155,-7){$\frac{\pi}{mk}$}
            \put(-142,-7){$\frac{2\pi}{mk}$} 
            \put(-0,-7){$2\pi$} \put(-168,153){$2\pi$}\put(-191,130){$\frac{(2ml-1)\pi}{ml}$}}
 \end{center}
    \caption{Dynamics of $Y_0$.}
    \label{figure:torus-general}
\end{figure}

This global result will be proved not only for system \eqref{equ100}, but more general $C^\infty$-smooth slow-fast systems on $\mathbb{T}^2$ as well. Instead of two critical \rev{manifolds} one of which is normally attracting and the other one is normally repelling, we can have $m$ normally attracting closed critical \rev{manifolds} and $m$ normally repelling closed critical \rev{manifolds} where $m\in\mathbb{N}$ is arbitrary and fixed. More precisely, let us look at the following generalization of system \eqref{equ100}:
\begin{equation}
\label{equ101}
Y_\epsilon: \begin{vf}
        \dot{x} &=&\sin{m(ly-kx)}   \\
        \dot{y} &=&\epsilon
    \end{vf}
\end{equation}
where $m\in\mathbb{N}$, $k\in\mathbb{N}$, $l$ is a non-negative integer and the pair $(k,l)$ is relatively prime. It is not difficult to see that $Y_0$ has $2m$ disjoint simple closed curves of singularities on $\mathbb{T}^2$\rev{, and their union define the critical manifold of $Y_0$. Moreover, each of such curves} 
wraps $k$ times vertically around $\mathbb{T}^2$ and $l$ times horizontally. This means that each closed curve crosses the horizontal interval $[0,2\pi[\times\left\{0\right\}$ exactly $k$ times and the vertical interval $\left\{0\right\}\times[0,2\pi[$ $l$ times (see Figure \ref{figure:torus-general}). Therefore, we have $m$ normally attracting \rev{critical manifolds} denoted by $C_-^1,\dots,C_-^m$ and $m$ normally repelling \rev{critical manifolds} denoted by $C_+^1,\dots,C_+^m$ (note that $k,m>0$).
We shall show that the slow-fast system \eqref{equ101} has exactly $m$ hyperbolically attracting 
limit cycles and $m$ hyperbolically repelling 
limit cycles for each $\epsilon>0$ and $\epsilon\sim 0$.

\begin{theorem}
\label{theorem-realization}
There exist $\epsilon_0>0$ and tubular neighborhoods $\mathcal{U}_\pm^i$ of $C_\pm^i$ inside $\mathbb{T}^2$, with $i=1,\dots,m$, such that system $Y_\epsilon$ produces exactly $1$ limit cycle in $\mathcal{U}_\pm^i$, denoted by $\mathcal{O}_\pm^i$, for each $\epsilon\in]0,\epsilon_0]$. Each limit cycle $\mathcal{O}_\pm^i$ is hyperbolically attracting ($-$) or repelling ($+$), turns around the unknotted torus $\mathbb{T}^2$ $k$ times vertically and $l$ times horizontally, and tends in Hausdorff sense to $C_\pm^i$ in the limit $\epsilon\to 0$. Moreover, fixing $\epsilon\in]0,\epsilon_0]$, for any $\tau\in\mathbb{T}^2$ such that $\tau\notin \mathcal{O}_\pm^i$, we have that $\alpha(\tau)$ is one of the repelling cycles $\mathcal{O}_+^i$ and $\omega(\tau)$ is one of the attracting cycles $\mathcal{O}_-^i$. 
\end{theorem}
Theorem \ref{theorem-realization} will follow from Theorem \ref{theorem-main} of Section \ref{section-statement} stated in a more general framework. Moreover, we further conjecture in Section \ref{sec:contact} that one can even consider the non-hyperbolic case where the critical manifolds may have \emph{regular nilpotent contact points of finite order}.  
\smallskip

\begin{remark}\leavevmode\label{remark-Poincare}
\begin{itemize}[leftmargin=*]
    \item For $k,l>0$ relatively prime, notice that the slow flow of $Y_\epsilon$, that is $\lim_{\epsilon\to0}\frac{1}{\epsilon}Y_\epsilon|_{\mathcal{O}^i_\pm}$, is a translation on $\T^2$. It then follows from Theorem \ref{theorem-realization} (and e.g. \cite{broer2010dynamical}) that the limit cycles of $Y_\epsilon$ have rotation number $\frac{l}{k}$ and that $Y_\epsilon$ has no dense orbits.
    \item Related to the previous point, let $J$ denote the horizontal interval $J=[0,2\pi[\times\left\{0\right\}$ and let $\Pi:J\to J$ be the Poincar\'e map induced by the flow of $Y_\epsilon$ for $\epsilon>0$ sufficiently small. It follows that $\Pi$ is well-defined, particularly since $\dot y=\epsilon$. This observation, in principle, would allow us to instead consider $\epsilon$-perturbations of translations over $\frac{l}{k}$ on $\T^1=\S^1$. We prefer to not consider this route because in our general result of Theorem \ref{theorem-main} this Poincaré map is not necessarily well-defined as we do not restrict the sign of the slow flow. For example, for a $(k,l)$-torus knot, a Poincaré map on the torus would need to be defined as the first return map \emph{after $k$ vertical rotations}. So, if the flow on the critical manifolds have distinct directions, there would be points on $J$ that return to $J$ before turning vertically along the torus.
\end{itemize}

\end{remark}

\begin{figure}
    \centering
    \includegraphics{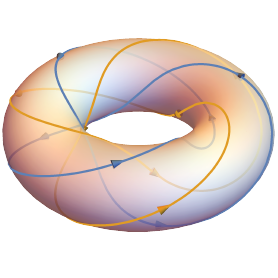}\hfill
    \includegraphics{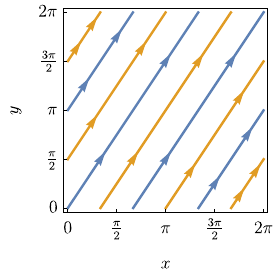}\\
    \includegraphics{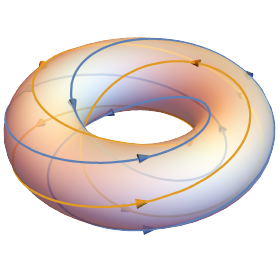}\hfill
    \includegraphics{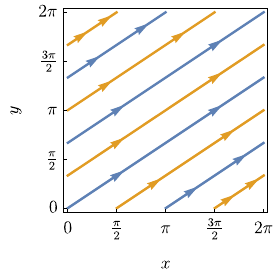}
    \includegraphics{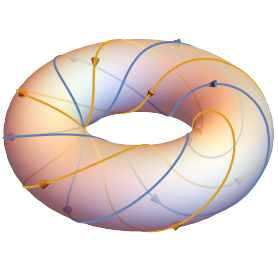}\hfill
    \includegraphics{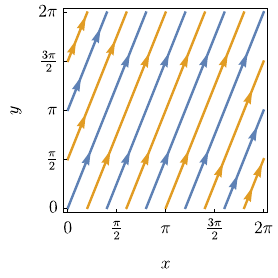}
    \includegraphics{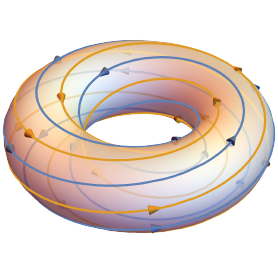}\hfill
    \includegraphics{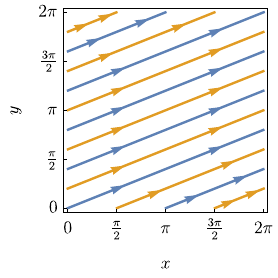}
    \caption{Examples of knotted critical manifolds  of \eqref{equ101} on $\T^2$ (left) and on $[0,2\pi]^2$ (right) for $m=1$: first row: $(3,2)$-knot (trefoil); second row: $(2,3)$-knot; third row $(5,2)$-knot (Solomon's seal); and fourth row: $(2,5)$-knot. Notice that $(k,l)$-knots are ambient isotopic (see Appendix \ref{sec:app1}) to $(l,k)$-knots.}
    \label{fig:crit_knots}
\end{figure}

A simple closed curve (i.e. an embedding $\mathbb{S}^1\to\mathbb{T}^2$) that turns around the torus $k$ times vertically (meridionally) and $l$ times horizontally (longitudinally) is often called a $(k,l)$-torus knot. We call a disjoint finite union of such torus knots a torus link (see Appendix \ref{sec:app1} and \cite{Adams,Knli,Rolf}). Using this terminology and Theorem \ref{theorem-realization}, we can say that for each small $\epsilon>0$ system $Y_\epsilon$ has a link consisting of $2m$ limit cycles (each limit cycle is a $(k,l)$-torus knot).  
 When $k=1$, $l=0$ or $l=1$, we deal with trivial knots or  unknotted circles (see, for example, system $X_\epsilon$). If $(k,l)=(3,2)$ or $(k,l)=(2,3)$, then the limit cycles of $Y_{\epsilon}$ are the trefoil knots (see Figure \ref{fig:crit_knots}). It is not possible to untangle the trefoil knot into the unknotted circle through $\mathbb{R}^3$-space without ``cutting" or ``gluing".  If $(k,l)=(5,2)$ or $(k,l)=(2,5)$, then $2m$ Solomon's seal knots occur in $Y_\epsilon$ for each $\epsilon>0$ (see Figure \ref{fig:crit_knots}). For more torus knots see a list in \cite{Rolf}.
 \smallskip



In \cite{FPR} the cyclicity of planar common slow-fast cycles has been studied. The common slow-fast cycles contain either attracting or repelling portions of the critical manifold. The authors were focused on the local study of a single common slow-fast cycle in the plane. In our paper we focus on the global study of slow-fast vector fields defined on $\mathbb{T}^2$ in the presence of a disjoint collection of common limit periodic sets of torus-knot type. 

\smallskip
In Section \ref{section-Fenichel} we define our slow-fast model on the $2$-torus and state the main result (Theorem \ref{theorem-main}). Section \ref{section-proofs} is devoted to the proof of Theorem \ref{theorem-main}. In Section \ref{sec:contact} we explain how \rev{one can} obtain a similar result to Theorem \ref{theorem-main} in the presence of regular nilpotent contact points.  In Section \ref{sec:app1} we give some basic definitions and results about torus knots.

\section{Assumptions and statement of results}
\label{section-Fenichel}

\subsection{Definition of a slow-fast model on $\mathbb{T}^2$ and assumptions}
Suppose that $X_{\epsilon,\rho}:\mathbb{T}^2\to T\mathbb{T}^2$ is a ($C^\infty$) smooth $(\epsilon,\rho)$-family of vector fields defined on the torus $\mathbb{T}^2$ of class $C^\infty$ where $\epsilon\ge 0$ is the singular parameter and $\rho\sim\rho_0\in\mathbb{R}^p$. The tangent bundle of $\mathbb{T}^2$ is denoted by $T\mathbb{T}^2$. The parameter $\rho$ is included for the sake of generality. The first assumption deals with the dynamics of the fast subsystem $X_{0,\rho}$.

\begin{assumption}
\label{assu-1}
The system $X_{0,\rho}$ has a smooth $\rho$-family of $2m$ disjoint embedded closed curves $C_{\rho,-}^1$, $C_{\rho,+}^1$, $\dots$, $C_{\rho,-}^m$, $C_{\rho,+}^m$ of singularities of $X_{0,\rho}$, for some $m\in\mathbb{N}$. Each curve wraps $k$ times vertically around the torus and $l$ times horizontally before it comes back to its initial point, for some fixed relatively prime non-negative integers $k$ and $l$ with $k+l>0$. Moreover, $C_{\rho,-}^i$ (resp. $C_{\rho,+}^i$) consists of normally attracting (resp. repelling) singularities for each $i=1,\dots,m$.
\end{assumption}
We call the critical manifold of $X_{0,\rho}$ (i.e., the disjoint collection of $2m$ $(k,l)$-torus knots from Assumption \ref{assu-1}) the $2m$-link and denote it by $C_\rho$.

\begin{remark}
Notice that for $m\geq 2$, the vector field $X_{0,\rho}$ is well-defined only if the critical manifolds $C_{\rho,\pm}^i$, $i=1,\ldots,m$, of $X_{0,\rho}$ alternate each other stability-wise along the fast fibers.
\end{remark}

\begin{remark}

It is worth noting that torus knots and slow-fast dynamics interact in a non-trivial way. To be more precise, while the two pairs of torus knots of Figure \ref{fig:ex_knots} are equivalent (up to homeomorphism and even ambient isotopy) on $\T^2$ \cite{Rolf}, they lead to nonequivalent slow-fast systems, as in Figure \ref{fig:ex_knots_crit}. In our main Theorem \ref{theorem-main} we restrict to normally hyperbolic critical manifolds. This has the advantage of making the proof more concise. However, as we argue in Section \ref{sec:contact}, it is possible to extend the results of Theorem \ref{theorem-main} to some cases where the critical manifold is not normally hyperbolic.

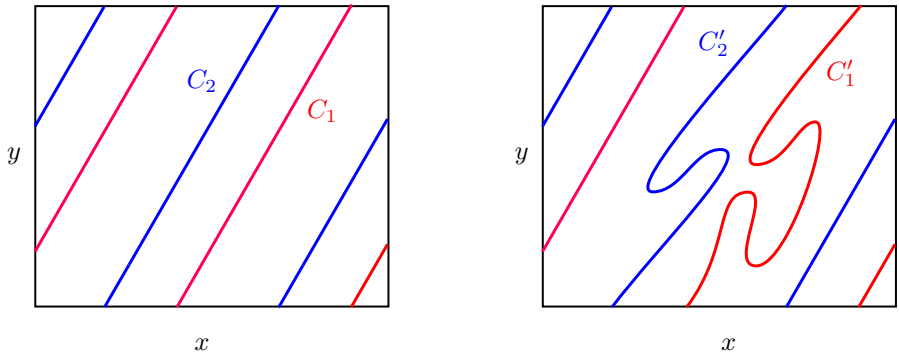
\begin{figure}[htbp]
    \centering
    \begin{tikzpicture}
    \node at (0,0) {
    \includegraphics[scale=0.75]{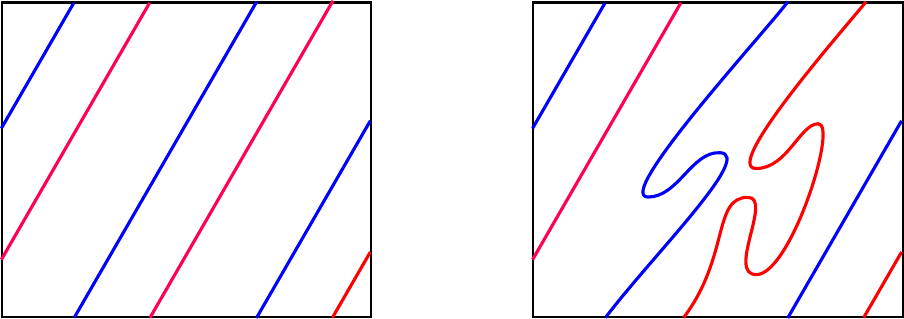}
    };
    \node at (-3.5,-2.5) {$x$};
    \node at (-6,0) {$y$};
    \node at ( 3.5,-2.5) {$x$};
    \node at (0.75,0) {$y$};
    \node at (-1.9,.6) {\color{red}{$C_1$}};
    \node at (-3.5,1) {\color{blue}{$C_2$}};
    \node at (5,1.1) {\color{red}{$C_1'$}};
    \node at (3.3,1.5) {\color{blue}{$C_2'$}};
    \end{tikzpicture}
    
    \caption{An example of a pair of equivalent (up to homeomorphism and even ambient isotopy) knots, that is, $C_i\sim C_i'$, $i=1,2$. Compare with Figure \ref{fig:ex_knots_crit}.}
    \label{fig:ex_knots}
\end{figure}

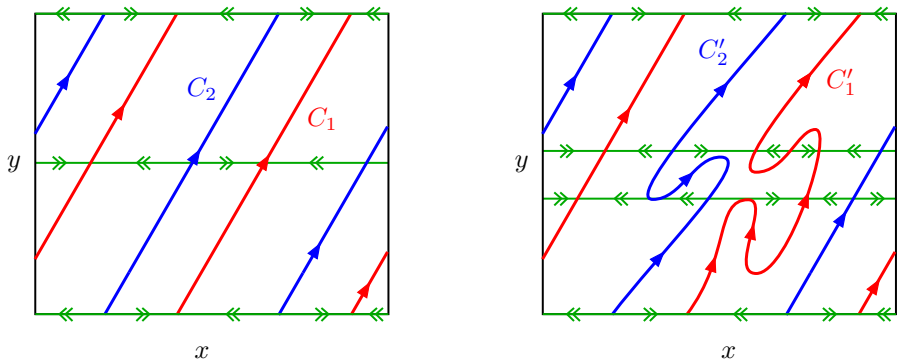
\begin{figure}[htbp]
    \centering
    \begin{tikzpicture}
    \node at (0,0) {
    \includegraphics[scale=0.75]{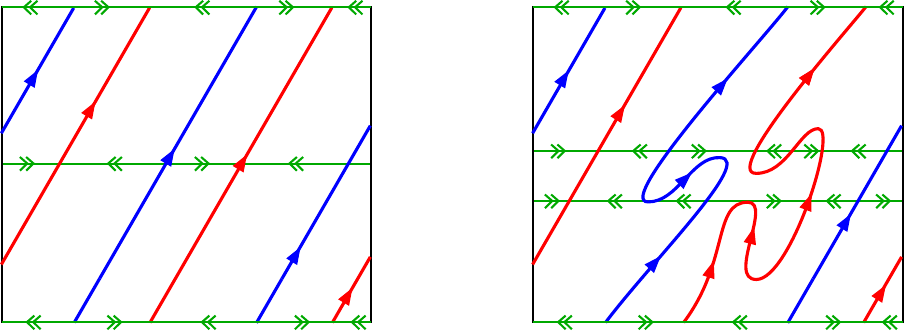}
    };
    \node at (-3.5,-2.5) {$x$};
    \node at (-6,0) {$y$};
    \node at ( 3.5,-2.5) {$x$};
    \node at (0.75,0) {$y$};
    \node at (-1.9,.6) {\color{red}{$C_1$}};
    \node at (-3.5,1) {\color{blue}{$C_2$}};
    \node at (5,1.1) {\color{red}{$C_1'$}};
    \node at (3.3,1.5) {\color{blue}{$C_2'$}};
    \end{tikzpicture}
    
    \caption{Even if the critical manifolds $C_i$ on the left and $C_i'$ on the right are equivalent as torus knots (compare with Figure \ref{fig:ex_knots}), they lead to completely different, and in fact nonequivalent, slow-fast dynamics. In the main part of this paper we assume that every knot of the $2m$-link critical manifold is normally hyperbolic, but refer also to Section \ref{sec:contact}.}
    \label{fig:ex_knots_crit}
\end{figure}


\end{remark}

A simple topological argument on $\mathbb{T}^2$ implies that disjoint torus knots have to be of the same type (see \cite{Rolf}). More precisely, we have
\begin{lemma}
If $C$ and $D$ are two disjoint knots on $\mathbb{T}^2$ of types $(k_1,l_1)$ and $(k_2,l_2)$ where $(k_1,l_1)$ and $(k_2,l_2)$ have the property given in Assumption \ref{assu-1}, then $k_1=k_2$ and $l_1=l_2$.
\end{lemma}

\rev{As we shall see shortly below, one important ingredient for our analysis is the \emph{slow-divergence integral}, defined in \eqref{eq:slow_di}.} In order to ensure that the slow divergence integral of $X_{\epsilon,\rho}$ is finite (i.e. well-defined), we have to assume that the slow dynamics of $X_{\epsilon,\rho}$ along the $2m$-link $C_\rho$ is regular. Let us recall that the slow dynamics or the slow flow, denoted by $\widetilde{X}_\rho$, is a $C^\infty$-smooth $\rho$-family of one-dimensional vector fields that describes the passage near $C_{\rho,-}^i$ or $C_{\rho,+}^i$, with $i=1,\dots,m$, when $\epsilon$ is positive and small. We have \[\widetilde{X}_\rho(\tau)=\lim_{\epsilon\to 0}\frac{(X_{\epsilon,\rho}+0\frac{\partial }{\partial \epsilon})(\tau)|_{M_\tau}}{\epsilon}\] where $M_\tau$ is any $C^n$ center manifold of $X_{\epsilon,\rho}+0\frac{\partial }{\partial \epsilon}$ at the point $\tau\in C_\rho$ ($n\in\mathbb{N}$ can be as large as we want). For more details about the definition of the slow dynamics see e.g. \cite{FPR}.


\begin{assumption}
\label{assu-2}
The slow dynamics $\widetilde{X}_{\rho_0}$ is nonzero at each point $\tau\in C_{\rho_0}$.
\end{assumption}
\begin{remark}
From Assumption \ref{assu-2} \rev{it} follows that along each of the $2m$ $(k,l)$-torus knots defined in Assumption \ref{assu-1} the slow dynamics $\widetilde{X}_{\rho_0}$ is regular (thus, without singularities). The slow dynamics may have different directions along the torus knots. For example, if we replace $\dot{y}=\epsilon$ in \eqref{equ100} with $\dot{y}=\epsilon \cos(y-x)$, then the slow flow along $C_-$ points upwards and along $C_+$ downwards. 
\end{remark}
\smallskip

\rev{Now that we have made the appropriate assumptions,} we can define the slow divergence integral of \rev{$X_{0,\rho}$} along $C_{\rho,-}^i$ (resp. $C_{\rho,+}^i$) as:
\begin{equation}\label{eq:slow_di}
    I_-^i(\rho)=\int_{C_{\rho,-}^i}\divergenceOperator X_{0,\rho}ds \quad \Big(\text{resp. }I_+^i(\rho)=\int_{C_{\rho,+}^i}\divergenceOperator X_{0,\rho}ds\Big)
\end{equation}
where $i=1,\dots,m$ and $s$ is the slow time of $\widetilde{X}_\rho$. 
    The slow divergence integral $I_{\pm}^i$ is independent of the local chart and the chosen volume form on $\mathbb{T}^2$, and represents the leading order term of the integral of divergence along orbits of  $X_{\epsilon,\rho}$, multiplied by $\epsilon$ (see Section \ref{section-proofs}). We typically compute $I_{-}^i$ by dividing the compact curve $C_{\rho,-}^i$ into a finite number of segments $[\tau_1,\tau_2],[\tau_2,\tau_3],\dots,[\tau_{r-1},\tau_{r}]$, where $r\in\mathbb{N}$, $\tau_1,\dots,\tau_r\in C_{\rho,-}^i$ and $\tau_1=\tau_r$ ($C_{\rho,-}^i$ is closed for all $\rho\sim\rho_0$), and by calculating the integral on each segment in suitable normal form coordinates:
    \begin{equation}\label{integral-division}
    I_-^i(\rho)=\int_{\tau_1}^{\tau_2}\divergenceOperator X_{0,\rho}ds+\int_{\tau_2}^{\tau_3}\divergenceOperator X_{0,\rho}ds+\cdots+\int_{\tau_{r-1}}^{\tau_r}\divergenceOperator X_{0,\rho}ds.
    \end{equation}
We assume that the points $\tau_1,\tau_2,\dots$ follow the direction of the slow flow along $C_{\rho,-}^i$ and that the above normally attracting segments are small enough such that on each segment we can use a Takens normal form for $C^n$-equivalence
\begin{equation}
\label{takens}
 \begin{vf}
        \dot{x} &=&-x   \\
        \dot{y} &=&\epsilon 
    \end{vf}
\end{equation}
(see \cite{DR}). The slow segment is lying inside $\{x=0\}$. Since the divergence of the vector field in \eqref{takens} for $\epsilon=0$ is $-1$, it is clear that each integral in \eqref{integral-division} is negative. Thus, the slow divergence integral $I_-^i$ is negative for all $\rho\sim\rho_0$ and $i=1,\dots,m$.
Similarly, we see that $I_+^i$ is positive, for all $\rho\sim\rho_0$ and $i=1,\dots,m$, because each $C_{\rho,+}^i$ is repelling.
\smallskip


\subsection{Statement of results}\label{section-statement}
In this section we state our main result. Let $X_{\epsilon,\rho}$ satisfy Assumptions \ref{assu-1}--\ref{assu-2}. Then for all $\epsilon\sim 0$ and $\epsilon>0$ a $2m$-link of type $(k,l)$ occurs inside $X_{\epsilon,\rho}$. More precisely,
\begin{theorem}
\label{theorem-main}
Suppose that system $X_{\epsilon,\rho}$ satisfies Assumptions \ref{assu-1}--\ref{assu-2}. Then there exist $\epsilon_0>0$, a neighborhood $\mathcal{V}$ of $\rho_0$ and tubular neighborhoods $\mathcal{U}_\pm^i$ of $C_{\rho_0,\pm}^i$, for each $i=1,\dots,m$, such that $X_{\epsilon,\rho}$ has exactly one limit cycle in $\mathcal{U}_\pm^i$, denoted by $\mathcal{O}_{\epsilon,\rho,\pm}^i$, for all $(\epsilon,\rho)\in ]0,\epsilon_0]\times\mathcal{V}$ and $i=1,\dots,m$. The limit cycle $\mathcal{O}_{\epsilon,\rho,-}^i$ (resp. $\mathcal{O}_{\epsilon,\rho,+}^i$) is a hyperbolic and attracting (resp. repelling) limit cycle. Moreover, each $\mathcal{O}_{\epsilon,\rho,\pm}^i$ tends in Hausdorff sense to the $(k,l)$-torus knot $C_{\rho_0,\pm}^i$ as $(\epsilon,\rho)\to (0,\rho_0)$. If $\tau$ is any point on $\mathbb{T}^2$ not lying in $\mathcal{O}_{\epsilon,\rho,\pm}^i$, then the $\omega$-limit (resp. the $\alpha$-limit) of $\tau$ is one of the attracting (resp. repelling) limit cycles $\mathcal{O}_{\epsilon,\rho,-}^1,\dots, \mathcal{O}_{\epsilon,\rho,-}^m$ (resp. $\mathcal{O}_{\epsilon,\rho,+}^1,\dots, \mathcal{O}_{\epsilon,\rho,+}^m$).
\end{theorem}
We prove Theorem \ref{theorem-main} in Section \ref{section-proofs}.

\section{Proof of Theorem \ref{theorem-main}}
\label{section-proofs}
This section is devoted to the proof of Theorem \ref{theorem-main}. Let $X_{\epsilon,\rho}$ satisfy Assumptions \ref{assu-1}--\ref{assu-2}. We divide the proof of Theorem \ref{theorem-main} into three parts:
\begin{figure}[htb]
    \begin{center}
        \ifnum\pdfoutput=0\includegraphics[bb=0 0 1302 839,width=11cm,height=5cm]{torusexample.png}\else\includegraphics[width=8.5cm,height=5.5cm]{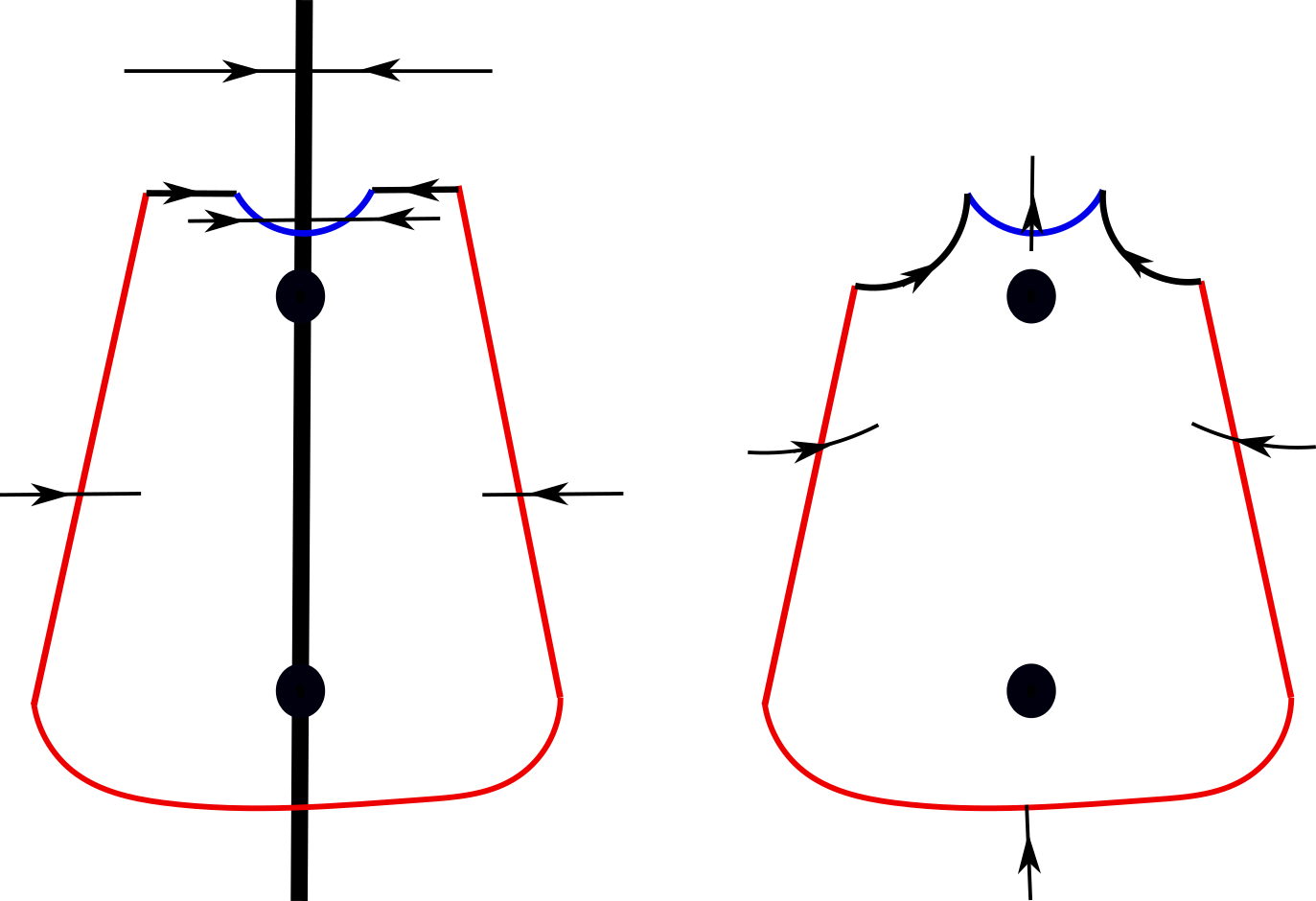}\fi
        {\footnotesize
            \put(-190,-8){$(a)$}  \put(-182,27){$\tau_j$}\put(-182,100){$\tau_{j+1}$}
            \put(-47,27){$\tau_j$}\put(-47,100){$\tau_{j+1}$}
            \put(-56,-8){$(b)$}  }
 \end{center}
    \caption{A flow box neighborhood with the inset (red) and the outset (blue).  (a) $\epsilon=0$ (b) $\epsilon>0$. \rev{Notice that for $\epsilon>0$ the vector field is transverse to the flow box. This essential property is lost when $\epsilon=0$.} }
    \label{figure:flowbox}
\end{figure}
\begin{enumerate}
    \item\textbf{(Flow-box neighborhoods near $C_{\rho_0,\pm}^i$)} Let $i\in\{1,\dots,m\}$. We construct a succession of flow box neighborhoods along the attracting \rev{critical manifold} $C_{\rho_0,-}^i$ (we can do the same with the repelling \rev{critical manifold} $C_{\rho_0,+}^i$ by reversing time). Since $X_{\epsilon,\rho}$ has no singularities on $\mathbb{T}^2$ for all small $\epsilon>0$ and $\rho\sim\rho_0$ (Assumptions \ref{assu-1}--\ref{assu-2} are satisfied), we can cover $C_{\rho_0,-}^i$ with a finite number of flow boxes that are uniform in $\epsilon>0$ and $\rho\sim\rho_0$, i.e. their size is fixed as $(\epsilon,\rho)\to (0,\rho_0)$  \rev{(see a sketch of a flow-box in Figure \ref{figure:flowbox}).} This, together with a fixed point theorem for segments, will enable us to show the existence of a closed orbit in an $(\epsilon,\rho)$-uniform tubular neighborhood  $\mathcal{U}_-^i$ of $C_{\rho_0,-}^i$, for each small $\epsilon>0$ and $\rho\sim\rho_0$. Our approach is based on the techniques from \cite{FPR}.
    
    \item \textbf{(The slow divergence integral along $C_{\rho_0,\pm}^i$)} We relate the slow divergence integral defined along $C_{\rho_0,-}^i$ (Section \ref{section-Fenichel}) to the integral of divergence along orbits of $X_{\epsilon,\rho}$ inside $\mathcal{U}_-^i$. Using a well-known connection between the derivative of the Poincar\'{e} map and the divergence integral, and the fact that the slow divergence integral is negative, we show that $X_{\epsilon,\rho}$ has at most one limit cycle in $\mathcal{U}_-^i$ up to shrinking $\mathcal{U}_-^i$ if needed. The limit cycle is a hyperbolically attracting $(k,l)$-torus knot. We also use a generalization of \rev{the} Poincar\'{e}-Bendixson Theorem to $\mathbb{T}^2$  due to Schwartz \cite{Schw}.
    
    \item \textbf{(Global dynamics of $X_{\epsilon,\rho}$ on $\mathbb{T}^2$)} To study the global dynamics on $\mathbb{T}^2$, we use the same result due to Schwartz \cite{Schw}. 
\end{enumerate}
\textit{1. Flow-box neighborhoods near $C_{\rho_0,\pm}^i$.} Let $\tau\subset\mathbb{T}^2$ and let $X$ be a $C^\infty$-smooth vector field on $\mathbb{T}^2$. We say that a neighborhood $\mathcal{U}$ of $\tau$ is a flow box neighborhood of $\tau$ for $X$ if $\mathcal{U}=F([0,1]\times [0,1])$ where $F:[0,1]\times[0,1]\to \mathbb{T}^2$ is a $C^\infty$-smooth diffeomorphism with the following property: the vector field $X$ is directed from outside to inside $\mathcal{U}$ along $F([0,1]\times\{0\})$ (the inset), from inside to outside $\mathcal{U}$ along $F([0,1]\times\{1\})$ (the outset), $F(\{0,1\}\times [0,1])$ are parts of orbits of $X$ and $X$ has no singularities in $\mathcal{U}$. Thus, all orbits of $X$ starting at the inset reach the outset. The set $\tau$ can be a point, a segment, etc.  We use this definition to introduce a slow-fast family of flow box neighborhoods of $\tau$ for $X_{\epsilon,\rho}$, i.e. a family $\{\mathcal{U}_{\epsilon,\rho};(\epsilon,\rho)\in [0,\epsilon_0]\times\mathcal{V}\}$ of neighborhoods of $\tau$ where $\epsilon_0>0$, $\mathcal{V}$ is a neighborhood of $\rho_0$, $\mathcal{U}_{\epsilon,\rho}=F_{\epsilon,\rho}([0,1]\times [0,1])$, $F_{\epsilon,\rho}:[0,1]\times[0,1]\to \mathbb{T}^2$ is an $(\epsilon,\rho)$-family of smooth diffeomorphisms such that $\mathcal{U}_{\epsilon,\rho}$ is a flow box neighborhood of $\tau$ for $X_{\epsilon,\rho}$ for each $(\epsilon,\rho)\in ]0,\epsilon_0]\times\mathcal{V}$ and such that the intersection of the sets $\mathcal{U}_{\epsilon,\rho}$, with $(\epsilon,\rho)\in [0,\epsilon_0]\times\mathcal{V}$, is a neighborhood of $\tau$. For more details see Definition 9 in \cite{FPR}.

\begin{figure}[htb]
    \begin{center}
        \ifnum\pdfoutput=0\includegraphics[bb=0 0 1302 839,width=11cm,height=5cm]{torusexample.png}\else\includegraphics[width=6cm,height=5.5cm]{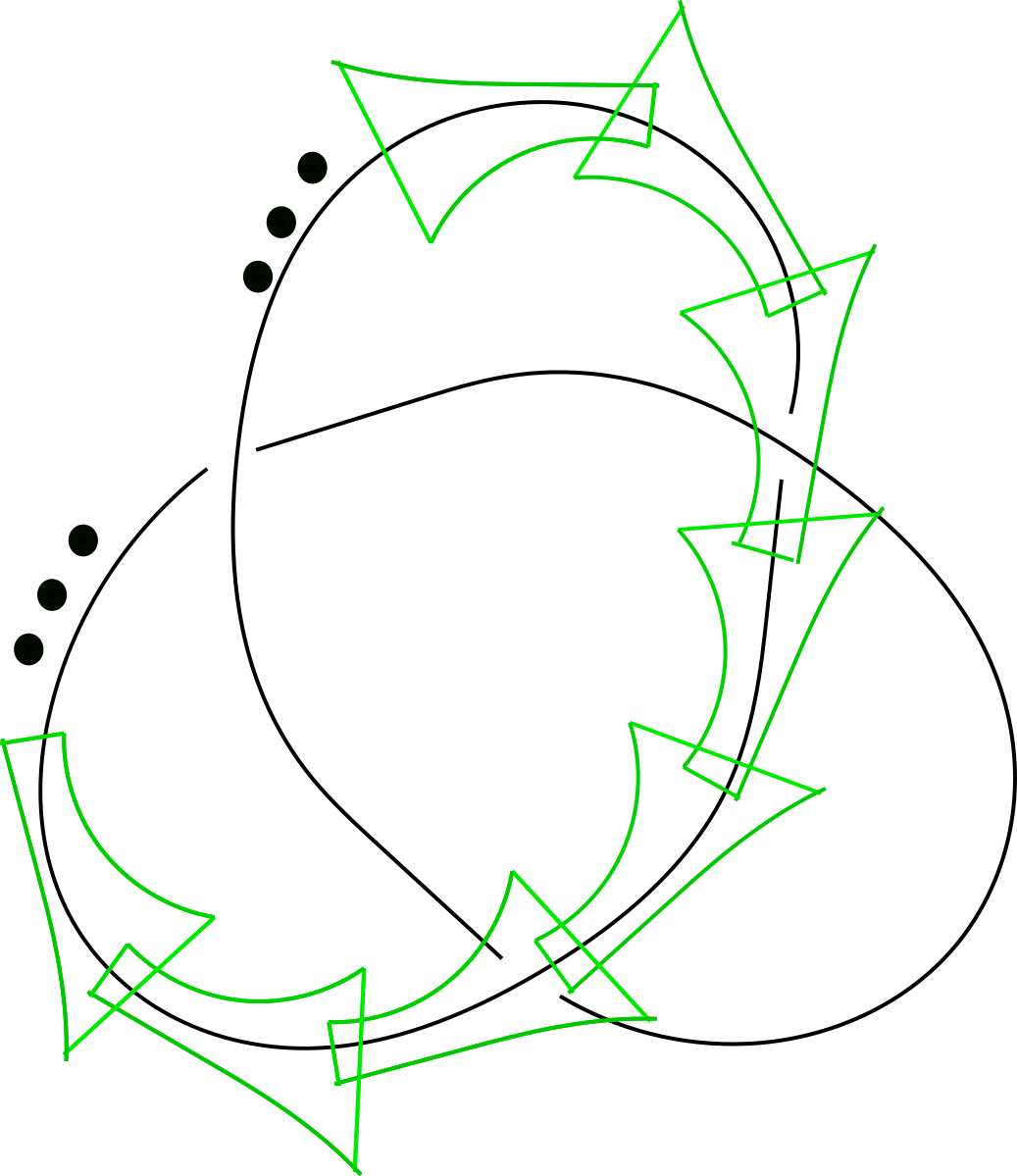}\fi
 \end{center}
    \caption{A trefoil knot covered with flow box neighborhoods.}
    \label{figure:sequence}
\end{figure}
Fix $i\in\{1,\dots,m\}$. We divide the compact \rev{critical manifold} $C_{\rho_0,-}^i$ into the segments $[\tau_1,\tau_2],[\tau_2,\tau_3],\dots,[\tau_{r-1},\tau_r]$ such that $\tau_j\in C_{\rho_0,-}^i$ for all $j=1,\dots,r$, $\tau_1=\tau_r$ and such that the points $\tau_1,\tau_2,\dots,\tau_r$ follow the direction of the slow flow along $C_{\rho_0,-}^i$ (i.e. the slow flow of $X_{0,\rho_0}$ goes from $\tau_1$ to $\tau_2$, from $\tau_2$ to $\tau_3$, etc.). For each segment $[\tau_j,\tau_{j+1}]\subset C_{\rho_0,-}^i$, $j=1,\dots,r-1$, we define a slow-fast family of flow box neighborhoods of $[\tau_j,\tau_{j+1}]$ inside $X_{\epsilon,\rho}$ (see Figure \ref{figure:flowbox}). We describe it using equivalence normal form coordinates. There exists a local chart on $\mathbb{T}^2$ around $[\tau_j,\tau_{j+1}]$ in which $X_{\epsilon,\rho}$ is locally given by the Takens normal form \eqref{takens} (the segment $[\tau_j,\tau_{j+1}]$ can be as small as we need). In the normal form coordinates, the segment is given by $[-1,1]$ on the $y$-axis. The outset is a parabola like segment above the line $\{y=1\}$. From the end points of the outset we have two parts of orbits of $X_{\epsilon,\rho}$ (they are horizontal in the limit $\epsilon=0$). The inset consists of two lines and a convex segment between them as indicated in Figure \ref{figure:flowbox}. Note that the vector field \eqref{takens} is directed from outside to inside the neighborhood along the inset for each small $\epsilon>0$ because the slow dynamics $y'=1$ points upwards near the $y$-axis. A detailed description of the slow-fast family of flow box neighborhoods can be found in Section 5 of \cite{FPR}.

Thus, we have constructed a finite number of flow box neighborhoods $\mathcal{U}_{\epsilon,\rho}^j$, $j=1,\dots,r-1$, along $C_{\rho_0,-}^i$ (see Figure \ref{figure:sequence}). Since each flow box neighborhood is contracted along its segment, we may assume that the outset of $\mathcal{U}_{\epsilon,\rho}^j$ is completely inside $\mathcal{U}_{\epsilon,\rho}^{j+1}$, for each $j=1,\dots,r-2$, and that the outset of $\mathcal{U}_{\epsilon,\rho}^{r-1}$ lies inside $\mathcal{U}_{\epsilon,\rho}^1$ (see Proposition 4 of \cite{FPR}). This implies that the Poincar\'{e} map from the inset of $\mathcal{U}_{\epsilon,\rho}^1$ to itself is well defined and smooth for $\epsilon>0$ small enough. Since the inset is a segment, there is a fixed point by Brouwer's Fixed Point Theorem. Thus, there exists a closed orbit near $C_{\rho_0,-}^i$ for each small $\epsilon>0$ and $\rho\sim\rho_0$.
\\
\\
\textit{2. The slow divergence integral along $C_{\rho_0,\pm}^i$.} Since $C_{\rho_0,-}^i$ is normally attracting and compact, it is well known (see e.g. \cite{DR}) that for any small $\kappa>0$ we can find $\epsilon_0>0$, a neighborhood $\mathcal{V}$ of $\rho_0$ and a tubular neighborhood $\mathcal{U}_-^i$ of $C_{\rho_0,-}^i$ such that
\begin{equation}
\label{estimate} 
\frac{I_-^i(\rho_0)-\kappa}{\epsilon}\le \int_{\mathcal{O}_{\epsilon,\rho}}\divergenceOperator{X_{\epsilon,\rho}}dt\le \frac{I_-^i(\rho_0)+\kappa}{\epsilon}
\end{equation}
where $I_-^i(\rho_0)<0$ is the slow divergence integral defined in Section \ref{section-Fenichel} and $\mathcal{O}_{\epsilon,\rho}$ is an arbitrary closed orbit of $X_{\epsilon,\rho}$ in $\mathcal{U}_-^i$ with $\epsilon\in]0,\epsilon_0]$ and $\rho\in\mathcal{V}$. We use now \eqref{estimate} to show that any closed orbit $\mathcal{O}_{\epsilon,\rho}$ of $X_{\epsilon,\rho}$ in $\mathcal{U}_-^i$ is hyperbolically attracting. It suffices
to take a $\kappa>0$ such that $\kappa<-I_-^i(\rho_0)$ and to use the Poincar\'{e} formula $P'(s_0)=e^{\int_{\mathcal{O}_{\epsilon,\rho}}\divergenceOperator{X_{\epsilon,\rho}}dt}$ where $P$ is the Poincar\'{e} map defined on a transverse section near $\mathcal{O}_{\epsilon,\rho}$ parametrized by a regular parameter $s$ ($\mathcal{O}_{\epsilon,\rho}$ corresponds to $s_0$). Note that the torus is orientable. Since the integral of divergence in \eqref{estimate} is negative (and of order $O(\frac{1}{\epsilon})$), we have $|P'(s_0)|<1$ and the closed orbit $\mathcal{O}_{\epsilon,\rho}$ is a hyperbolically attracting limit cycle.

To prove that $C_{\rho_0,-}^i$ can produce at most one limit cycle inside $\mathcal{U}_-^i$, we use the following result due to Schwartz \cite{Schw} (see also Theorem 6.6 in \cite{Cie}):
\begin{theorem}
\label{theorem-sch}
Assume that $\dot{x}=X(x)$ is an autonomous system of class $C^2$ on a
compact, connected two-dimensional orientable manifold $M$ of class $C^2$. Assume that this system defines a flow and that the manifold $M$
is not a minimal set. Then, if the $\omega$-limit set $\omega(\tau)$ of a point $\tau\in M$ does not contain
any critical point, $\omega(\tau)$ must be homeomorphic to the circle $\mathbb{S}^1$ (i.e. is a periodic trajectory).
\end{theorem}
We apply Theorem \ref{theorem-sch} to $M=\mathbb{T}^2$ and $X=-X_{\epsilon,\rho}$ for each fixed $\epsilon>0$ and $\rho\sim\rho_0$ (in this paper, both $\mathbb{T}^2$ and $X_{\epsilon,\rho}$ are of class $C^\infty$). Let us recall that a set $S\subset M$ is minimal if it is nonempty, invariant, closed and there is no proper subset of $S$ that has all these properties. Since $X_{\epsilon,\rho}$ has at least one closed (or periodic) orbit for each $\epsilon>0$ small enough (see Step 1), $\mathbb{T}^2$ is not minimal. As already mentioned, $X_{\epsilon,\rho}$ has no critical points on $\mathbb{T}^2$ for $\rho\sim\rho_0$ and $\epsilon>0$ small enough.

Suppose now that $X_{\epsilon,\rho}$ has at least two limit cycles in $\mathcal{U}_-^i$ ($\mathcal{O}_1$ and $\mathcal{O}_2$). Then $\mathcal{O}_1$ and $\mathcal{O}_2$ are $(k,l)$-torus knots and bound an invariant set $S$. Take any $\tau\in S$ not lying on a closed orbit of $X_{\epsilon,\rho}$ (such $\tau$ exists if we are close enough to the isolated closed orbits $\mathcal{O}_1$ or $\mathcal{O}_2$). Following Theorem \ref{theorem-sch}, $\omega(\tau)\subset S$ w.r.t. $-X_{\epsilon,\rho}$ is a periodic trajectory different from $\mathcal{O}_1$ and $\mathcal{O}_2$ ($\mathcal{O}_1$ and $\mathcal{O}_2$ are repelling for $-X_{\epsilon,\rho}$). This leads to a contradiction because the periodic trajectory $\omega(\tau)$ is not repelling (all closed orbits generated by $C_{\rho_0,-}^i$ are hyperbolic and repelling limit cycles w.r.t. $-X_{\epsilon,\rho}$). 
\smallskip

From Step 1 and Step 2 it follows that there exists $\epsilon_0>0$, a small neighborhood $\mathcal{V}$ of $\rho_0$ and a tubular neighborhood $\mathcal{U}_-^i$ of $C_{\rho_0,-}^i$ such that for each $(\epsilon,\rho)\in ]0,\epsilon_0]\times  \mathcal{V}$ system $X_{\epsilon,\rho}$ has precisely one limit cycle in $\mathcal{U}_-^i$. It is hyperbolic and attracting and clearly of $(k,l)$-torus knot type. \rev{Moreover, using the Takens normal form \eqref{takens} and exponential contraction of orbits near $C_{\rho_0,-}^i$ as $\epsilon\to 0$, the limit cycle tends in Hausdorff sense to $C_{\rho_0,-}^i$ as $(\epsilon,\rho)\to (0,\rho_0)$ (see also Section 7.1 in \cite{de2021canard}).} By reversing time, we can prove a similar result near $C_{\rho_0,+}^i$ for each $i=1,\dots,m$ ($C_{\rho_0,+}^i$ generates one hyperbolic and repelling limit cycle). 
\\
\\
\textit{3. Global dynamics of $X_{\epsilon,\rho}$ on $\mathbb{T}^2$.} First note that the $2m$ limit cycles obtained in Step 1 and Step 2 are the only possible periodic trajectories of $X_{\epsilon,\rho}$. Indeed, any orbit with the initial point $\tau\in \mathbb{T}^2\setminus \cup_{i=1}^m\mathcal{U}_\pm^i$ ($\tau$ is thus uniformly away from the critical manifold of $X_{0,\rho_0}$) cannot be periodic (it is attracted to an attracting closed part $C_{\rho_0,-}^i$ and stays there forever). Theorem \ref{theorem-sch} implies now the rest of Theorem \ref{theorem-main}.

\begin{remark}
Naturally, there are other alternative arguments to show uniqueness of the limit cycle produced by $C_{\rho_0,-}^i$. One possibility is to use the Takens normal form \eqref{takens} and study the composition of exponentially strong linear contractions and $C^k$-smooth coordinate changes to show the uniqueness of the limit cycle, see for example \cite[Theorem 3]{GuIl}.
\end{remark}


\section{Regular nilpotent contact points of finite order}\label{sec:contact}





In this section we briefly argue that a similar result to Theorem \ref{theorem-main} can be obtained in some cases where the critical manifold $C_\rho$ is not normally hyperbolic. Instead of providing all the technicalities, we refer to \cite{FPR}, and only indicate the main ideas. \rev{As we do not present a rigorous proof, we acknowledge that these claims are a conjecture in the current work.}

\begin{definition}\label{def:contact}
A point $\tau\in C_{\rho_0}$ is called a \emph{nilpotent contact point} if the linear part of $X_{0,\rho_0}$ at $\tau$ (computed in local coordinates) is nilpotent. Near such a point, with $(\epsilon,\rho)\sim (0,\rho_0)$, system $X_{\epsilon,\rho}$ is $C^\infty$-equivalent to $\{\dot{x}=y-f(x,\rho),\dot{y}=\epsilon(g(x,\epsilon,\rho)+O(y-f(x,\rho))\}$ where $(x,y)$ are local normal form coordinates, $\tau$ is given by $(x,y)=(0,0)$, $f,g,O$ are $C^\infty$-smooth and the $1$-jet $j_1f(x,\rho_0)$ at $x=0$ is zero (see \cite{FPR}).  A nilpotent contact point $\tau$ is further called \emph{regular} if $g(0,0,\rho_0)\ne 0$. Finally, a nilpotent contact point $\tau$ is called of finite order, if there exists an integer $n\geq2$ such that the vanishing order of $f(x,\rho_0)$ at $x=0$ is equal to $n$. 
\end{definition}
It is not difficult to see that the above definition is independent of the chosen normal form. \rev{We focus on limit cycles generated by so-called regular common cycles (they contain either only attracting or only repelling portions of the critical manifold, all involved
contact points are regular and the slow dynamics on all portions has no singular points, see \cite[Definition 4.4]{de2021canard}). A regular
common cycle can contain different elementary slow–fast segments as explained in \cite[Figure 4.9]{de2021canard} or \cite{FPR}. For example, near a ``regular nilpotent point of \emph{even} contact order'' we can have a jump point (see Figure \ref{fig:example2}). Notice that canard points do not fall into Definition \ref{def:contact}.}

It turns out that the main ideas of Theorem \ref{theorem-main} also hold\rev{, under appropriate assumptions, } to the case where the critical manifold may contain isolated regular nilpotent contact points of finite order. To be more precise, Assumption 1 becomes:

\begin{description}[before={\renewcommand\makelabel[1]{\bfseries ##1.}},leftmargin=*]
\item[Assumption 1']\textit{
The system $X_{0,\rho}$ has a smooth $\rho$-family of $2m$ disjoint embedded curves $C_{\rho,-}^1$, $C_{\rho,+}^1$, $\dots$, $C_{\rho,-}^m$, $C_{\rho,+}^m$ of singularities of $X_{0,\rho}$, for some $m\in\mathbb{N}$. Each curve wraps $k$ times vertically around the torus and $l$ times horizontally before it comes back to its initial point, for some fixed relatively prime positive integers $k$ and $l$. Moreover, $C_{\rho,\pm}^i$ consists of normally hyperbolic singularities for each $i=1,\dots,m$ except at a finite number of points where $C_{\rho,-}^i$ (resp. $C_{\rho,+}^i$) may have isolated regular nilpotent contact points of finite order. \rev{Furthermore, we assume that if a singular jump occurs (either in forward or backward time), then it occurs within the same critical manifold $C^i_{\rho,\pm}$. In other words, we do not allow configurations of the critical manifold where a contact point of $C^i_{\rho,\pm}$ connects via the fast dynamics to a different (disjoint) component of the critical manifold.  }
}
\end{description}

\begin{figure}[htbp]
    \centering
    \begin{tikzpicture}
    \node at (-3,0){
    \includegraphics{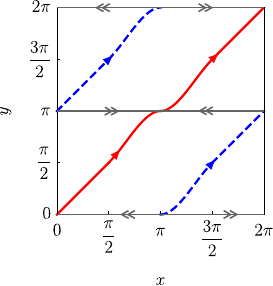}
    };
    \node at (3,0){
    \includegraphics{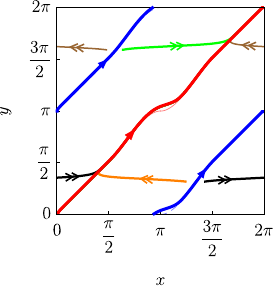}
    };
    \end{tikzpicture}
    
    \caption{ Example of a $2$-link critical manifold consisting of two $(1,1)$-knots with an odd regular contact point. An odd contact point is locally given by the singularity $y=x^{2n+1}$, $n\in\mathbb N$. On the left we sketch the singular limit, observe that the critical manifold  has an odd contact with the fast foliation at $x=\pi$. The red knot corresponds to the attracting critical manifold, while the blue knot to the repelling one. On the right we sketch possible orbits for $\epsilon>0$ sufficiently small, where the red knot indicates the attracting  limit cycle and the  blue knot the repelling limit cycle. Orbits away from the \rev{repelling} limit cycle approach the \rev{attracting} limit cycle as $t\to\infty$.}
    \label{fig:example1}
\end{figure}

\begin{figure}[htbp]
    \centering
    \begin{tikzpicture}
    \node at (-3,0){
    \includegraphics{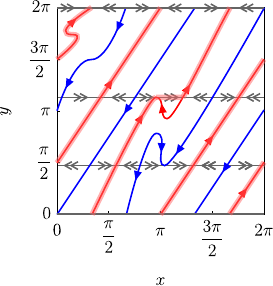}
    };
    \node at (3,0){
    \includegraphics{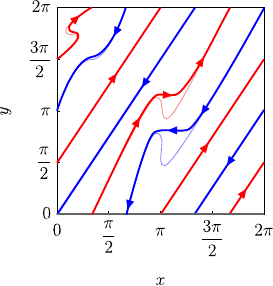}
    };
    \end{tikzpicture}
    
    \caption{ Example of a $2$-link critical manifold consisting of two Trefoils ($(3,2)$-knots) with regular contact points. On the left we sketch the singular limit, observe that the critical manifolds (thick red and blue knots) contain, for example, jump points. All flows are shown in forward time. The shaded red path indicates a singular slow-fast cycle, that is a knot consisting of slow and fast orbits of the singular limit. As sketched on the right, the singular slow-fast knot perturbs to an attracting limit cycle, which in fact corresponds to a relaxation oscillation. On the other hand, the blue knot perturbs to a repelling limit cycle, also of slow-fast type (singular limit not shown). The main argument to show that a unique limit cycle bifurcates from each singular slow-fast knot for $\epsilon>0$ sufficiently small, is that the local passage through the regular contact points is a contraction.}
    \label{fig:example2}
\end{figure}

Then, one can \rev{conjecture} that under Assumptions 1' and 2 (Assumption 2 is valid at normally hyperbolic singularities and the slow flow is unbounded at the contact points), $X_{\epsilon,\rho}$ has $2m$ limit cycles \rev{of $(k,l)$-torus knot type} from which $m$ are attracting and $m$ are repelling. We believe that the proof would follow similar arguments as in Section \ref{section-proofs} and \cite{FPR}. 
Examples of knotted critical manifolds with regular nilpotent contact points of finite order are provided in Figures \ref{fig:example1} and \ref{fig:example2}.

Although we do not provide a specific model of a smooth slow-fast system on $\T^2$ with knotted critical manifold with contact points, we point-out that it is possible to construct such systems using, for example, Hermite interpolation \cite{spitzbart1960generalization}. 

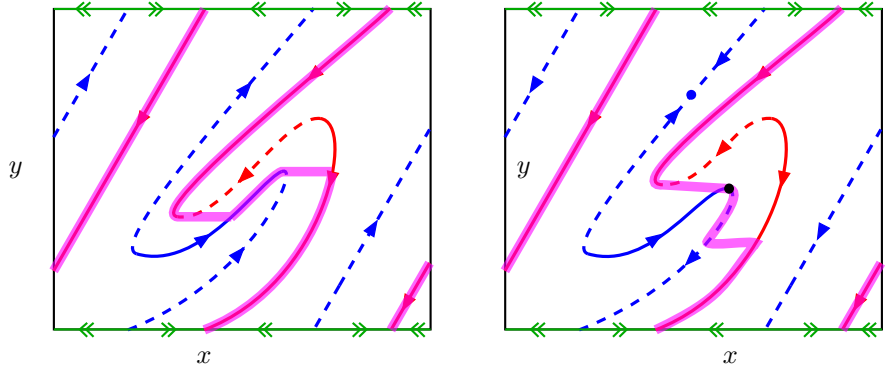
\begin{figure}[htbp]
    \centering
    \begin{tikzpicture}
    \node at (-3,0){
    \includegraphics[scale=0.8]{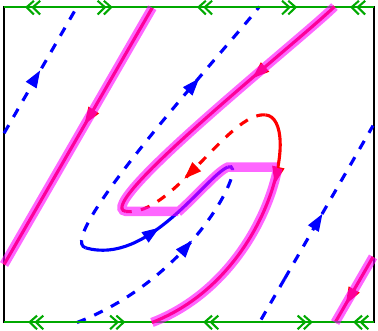}
    };
    \node at (3,0){
    \includegraphics[scale=0.8]{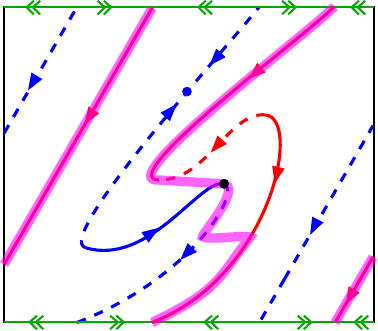}
    };
    \node at (-3.5,-2.5) {$x$};
    \node at (-6,0) {$y$};
    \node at ( 3.5,-2.5) {$x$};
    \node at (0.75,0) {$y$};
    \end{tikzpicture}
    \caption{\rev{We present in this figure two  singular closed candidate orbits (in magenta) not covered in our conjecture. On the left we sketch a slow-fast system on $\mathbb T^2$ with two critical manifolds having only jump points. Repelling regions are depicted as dashed curves, while attracting regions are solid. Notice that due to the configuration of the critical manifolds, it is possible that an orbit jumps from one critical manifold to another. In this particular case, we observe that the conjecture seems to hold, as from the singular limit one can clearly see that for $\epsilon>0$ sufficiently small, two limit cycles arise, the attracting one being a perturbation of the magenta candidate orbit. On the right we sketch yet another situation where the black dot is a canard point, thus not a regular nilpotent point, while all other contact points are regular. Notice that such candidate orbit would correspond to a canard cycle. The blue dot on the repelling blue critical manifold is introduced so that the slow flow is compatible with the canard point. Naturally, its position is arbitrary. We notice though that this canard situation can also arise without a canard point. These and other complicated situations are topics of future research.}}
    \label{fig:examples}
\end{figure}

\rev{
Finally, let us briefly notice the importance of Assumption 1': allowing jumps between different components of the critical manifold can lead to situations as depicted in Figure \ref{fig:examples} and which can definitely modify the number of limit cycles for $\epsilon>0$ sufficiently small. Furthermore, such a connection does not necessarily require canard points, see also Figure  \ref{fig:examples}, although these could introduce extra difficulties in the analysis.
}

\rev{
\section*{Acknowledgments}
The authors thank the anonymous reviewer for the constructive feedback provided, which helped to improve our manuscript.
}

\newpage 
\appendix

\section{Basic notions of Knot Theory}\label{sec:app1}




In this section we present some basic notions, definitions, and results about knot theory and, in particular, torus knots. For more details see, for example, \cite{Adams,Rolf}. We start with the definition of a knot:

\begin{definition} A \emph{knot} is an embedding $K:\S^1\to\R^3$ of a $1$-sphere into $\R^3$. A \emph{link} $L:\cup\, \S^1\to\R^3$ is a disjoint, finite union of knots\footnote{More generally, one can define a knot as being an embedding $K:\S^r\to X$, where $X$ is a topological space (usually $\S^n$ or $\R^n$). Alternatively, one can consider a knot as being a subset homemorphic to $\S^r$, $r\geq1$.}.
\end{definition}

In order to introduce a notion of equivalence between knots, we have the following:

\begin{definition}
A homotopy $h_t:\R^3\to\R^3$ is called an \emph{ambient isotopy} if $h_0$ is the identity and every $h_t$ is a homeomorphism.
\end{definition}

There are several notions of equivalence for knots. For our purposes it will suffice to keep in mind the following two:
\begin{definition}\leavevmode
\begin{enumerate}
    \item Two knots (or links) $K$ and $K'$ are \emph{equivalent}, if there is a homeomorphism $h:\R^3\to\R^3$ such that $h(K)=K'$.
    \item Two knots (or links) $K$ and $K'$ are \emph{ambient isotopic}, if there is an ambient isotopy such that $h_1(K)=K'$.
\end{enumerate}
\end{definition}

\begin{figure}[htbp]
    \centering
    \begin{tikzpicture}
    \node[] at (0,0){
    \includegraphics{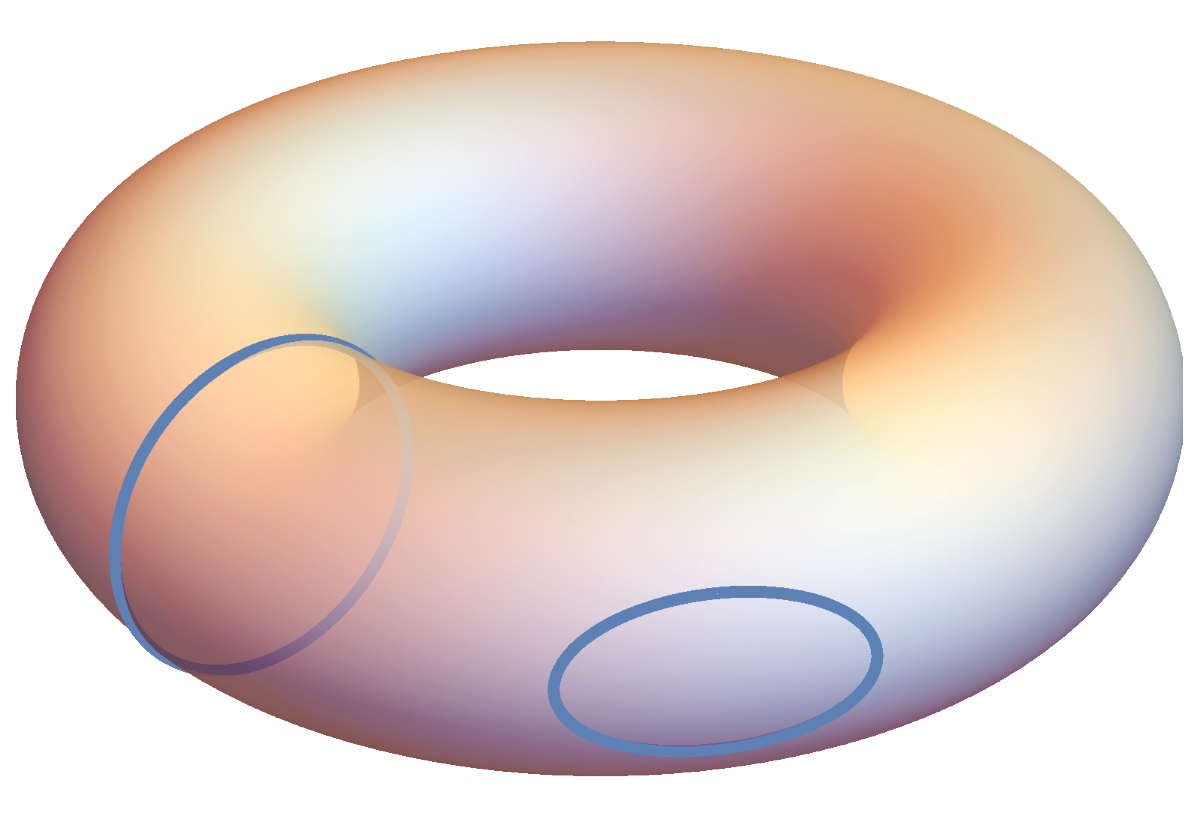}
    };
    \node at (-2.2,0){$K_1$};
    \node at (1,-.6){$K_2$};
    \end{tikzpicture}
    
    \caption{Representatives of the two classes of nonequivalent torus knots (up to homeomorphism). That is, every knot on the $2$-torus can be transformed via a homeomorphism to either $K_1$, which is a knot of class $(1,0)$, or $K_2$, which is a knot of class $(0,0)$ \cite{Rolf}.}
    \label{fig:equiv_knots}
\end{figure}

In this paper we are concerned with\textbf{} torus knots, that is, knots that lie in the $2$-torus $\T^2=\S^1\times\S^1\subset\R^3$. In this regard, we say that a torus knot $K$ is of class $(k,l)$ if the knot winds the torus $k$-times vertically and $l$-times horizontally, and where $k$ and $l$ are coprime (see e.g. Figure \ref{fig:crit_knots}). It turns out that the only nonequivalent knots (up to homeomorphism) are those presented in Figure \ref{fig:equiv_knots}.

We see that equivalence up to homeomorphism is too coarse. In contrast, ambient isotopies preserve orientation, which leads to a finer classification. 
\begin{theorem}[\cite{Knli,Rolf}]
Let $K_i$, $i=1,2$, be torus knots of class $(k_i,l_i)$ respectively. Then $K_1$ and $K_2$ are ambient isotopic if and only if $(k_1,l_1)=\pm(k_2,l_2)$. If one considers only the case where $k_i$ and $l_i$ are positive, then $K_1$ and $K_2$ are ambient isotopic if and only if $k_1=k_2$ and $l_1=l_2$ or $k_1=l_2$ and $l_1=k_2$.
\end{theorem}

Examples of equivalent knots modulo ambient isotopy can be  seen in Figure \ref{fig:crit_knots}.


\bibliographystyle{plain}
\bibliography{bibtex}
\end{document}